\documentstyle[]{article}
\title{Confidence Regions for Means of Multivariate Normal Distributions
and a Non-Symmetric Correlation Inequality for Gaussian Measure }
\author{Stanislaw J. Szarek \and Elisabeth Werner}
\date{}

\input amssym.def
\input amssym.tex
\input BoxedEPS
\SetepsfEPSFSpecial

\begin{document}
\maketitle
\begin{abstract}
\noindent Let $\mu$ be a Gaussian measure (say, on ${\bf R}^n$) and let $K,
L \subseteq {\bf R}^n$
be such that K is convex, $L$ is a ``layer"  (i.e.
 $L = \{ x : a \leq \langle x,u \rangle \leq b \}$ for some $a$, $b \in
{\bf R}$
and $u \in {\bf R}^n$) and the centers of mass (with respect to $\mu$) of
$K$ and $L$
coincide.  Then $\mu(K \cap L) \geq \mu(K) \cdot \mu(L)$.  This is
motivated by the
well-known  ``positive correlation conjecture" for symmetric sets and a related
inequality of Sidak concerning confidence regions for means of multivariate
normal
distributions.  The proof uses an apparently hitherto unknown estimate for the
(standard) Gaussian cumulative distribution function:
$\Phi (x) > 1 - \frac{(8/\pi)^{\frac{1}{2}}}{3x + (x^2 +8)^{\frac{1}{2}}}
e^{-x^2/2}$
(valid for $x > -1$).

\end{abstract}

\section{Introduction}

Let $\mu=\mu_n$ be the standard Gaussian measure on ${\bf R}^n$ with density
$(2 \pi )^{- \frac{n}{2}} e^{- \frac{\| x \|^2}{2}}$  (or any{\em centered}
Gaussian measure on ${\bf R}^n$).  It is a well known open problem
whether any two symmetric (with respect to the origin) convex sets $K_1$
and $K_2$ in ${\bf R}^n$ are positively
correlated
with respect to $\mu$, i.e. whether the following inequality holds
\begin{equation}
\mu (K_1 \cap K_2 ) \geq \mu (K_1) \mu (K_2).
\end{equation}
Of course once (1) is proved, it follows by induction that the following
formally stronger
statement is true:
\begin{equation}
\mu (K_1 \cap K_2 \cap \ldots \cap K_N ) \geq \mu (K_1) \mu (K_2)
\ldots \mu (K_N)
\end{equation}
for any convex symmetric sets  $K_1, K_2, \ldots, K_N$ in ${\bf R}^n$
(the same remark applies to any class of sets closed
under intersections).  In the language of statistics,  (1) and (2)
can be viewed as statements about
confidence regions for means of multivariate normal distributions
(cf. Theorem 1A below).
In some special cases (1) and  (2) are known to be true.  Pitt [P] proved
in 1977 that
(1) (hence (2)) holds in ${\bf R}^2$.  Also, if $K_1, K_2 \ldots,K_N$ are
symmetric layers in ${\bf R}^n$, i.e., sets of the form
$$K_i = \{ x \in {\bf R}^n : | \langle x, u_i \rangle | \leq 1 \}, u_i \in
{\bf R}^n, i =1,2,\ldots,N,$$
then (2)  holds  (note that in that particular case  (1)  {\em doesn't}
imply  (2)).
This was proved by Sidak [S] in 1967 and consequently is referred to as
Sidak's Lemma.  See also Gluskin [G] for a proof of Sidak's Lemma.
The proof gives in fact a version of (1) with $K_1$ - an arbitrary symmetric
convex body and  $K_2$ - a layer;  (2) {\em for layers} follows then by
induction.  We show in Remark 6 of Section 3 how
Sidak's Lemma can be proved easily with the approach of this paper
(an argument of this type seems to have recently occurred more or less
simultaneously
to several people).  In 1981 Borell [B] proved that (1) holds
for a class of convex symmetric bodies in ${\bf R}^n$ with certain
additional properties.  Recently Hu [H] proved a correlation inequality for
Gaussian measure involving convex functions rather than sets.
See also [S-S-Z] for a historical survey and other partial results and
[K-MS] for related results.

Here we prove the following.

\medskip\noindent {\bf Theorem 1.} {\em Let $K \subseteq {\bf R}^n$ be a
convex body and
$u \in {\bf R}^n \setminus \{0\}$ be such that
$$\int_K (\langle x,u \rangle - c) d\mu_n (x) =0$$
i.e. the centroid of $K$ with respect to $\mu_n$ lies on the hyperplane
$$H_c = \{ x \in  {\bf R}^n : \langle u,x \rangle =c \}.$$
\noindent Let $L = L(a,b) = \{ x \in {\bf R}^n : a \leq \langle x,u \rangle
\leq b \}$
where $a$, $b
\in {\bf R}$ are such that the centroid of $L$ also lies in $H_c$.  Then}
$$\mu_n (K \cap L) \geq \mu_n (K) \cdot \mu_n (L).$$

It is clear that Theorem 1 formally implies an analogous statement with $\mu_n$
replaced by {\em any} gaussian measure on ${\bf R}^n$ (centered or not).
In the language of
``confidence regions",  Theorem 1 may be restated as:

\medskip\noindent {\bf Theorem 1A.}  {\em Let} $X_1, X_2, \ldots, X_N, Y$
{\em be jointly Gaussian
random variables and} $b_1, b_2, \ldots, b_N, a, b \in {\bf R}$ {\em be
such that}
$${\bf E}(Y | X_1 \leq b_1, X_2 \leq b_2, \ldots, X_N \leq b_N, a \leq Y
\leq b) =
{\bf E}(Y | a \leq Y \leq b).$$
{\em Then}
$${\bf P}(X_1 \leq b_1, X_2 \leq b_2, \ldots, X_N \leq b_N, a \leq Y \leq b)$$
$$\geq \; \; {\bf P}(X_1 \leq b_1, X_2 \leq b_2, \ldots, X_N \leq b_N)
\cdot {\bf P}(a \leq Y \leq b).$$

\bigskip We point out that the discrepancy between the degrees of generality
 of Theorems 1 and 1A (general convex sets vs. ``rectangles")  is only apparent:
passing from rectangular to general parallelepipeds requires only a change
of variables;
a general convex polytope is a ``degenerated" parallelepiped,  and any
convex set can be
approximated by polytopes.

\medskip  Theorem 1  leads naturally to the following generalization of
the ``correlation conjecture" (1).

\medskip\noindent {\bf Problem 2.}  {\em If}  $K_1,
K_2 \in {\bf R}^n$ {\em are  convex sets  (not necessarily symmetric) such that
their centroids with respect to} $\mu_n$ {\em coincide,  does} (1) {\em hold?}

\medskip It is conceivable that the ``equality of the centroids" hypothesis
is not the most proper here and that it should be modified.  However,  we
were led
to that particular hypothesis while considering some variational arguments
related
to the original (symmetric) correlation conjecture (those arguments yield,
in particular, an alternative proof of the two-dimensional case shown in [P]).
Our Theorem 1 is related to Problem 2 in roughly the same way as Sidak's
Lemma is to
the original ``symmetric" conjecture.

\medskip Theorem 1  is proved in Section 3  (with proofs of some technical
lemmas
relegated to Section 4).  In Section 2 we develop some of the tools
necessary for
the proof.  They
may also be of independent interest, in particular Proposition 3 which
gives an upper estimate
on the tail of the Gaussian distribution that is sharper than the corresponding
``Komatsu inequality'' known from the literature (cf. [I-MK], p. 17;  see
also [Ba]
for another type of estimate).

\medskip \noindent {\bf Proposition 3.}  {\em  For} $x > -1$
$$\frac{2}{x + (x^2 +4)^{\frac{1}{2}}} \leq e^{\frac{x^2}{2}}
\int_x^{\infty} e^{- \frac{t^2}{2}}
dt \leq \frac{4}{3x + (x^2 +8)^{\frac{1}{2}}}$$

\medskip The lower estimate in Proposition 3 is the other ``Komatsu inequality"
and is true for any $x \in {\bf R}$.  The comparison of the upper estimate from
Proposition 3 with classical estimates is given in a table in Remark 4 in the
next section.

\medskip {\em Acknowledgement.} Research partially supported by authors'
respective grants from the National Science Foundation.  The final part of the
research has been performed while the authors were in
residence at MSRI Berkeley.  They express their gratitude to the staff of
the institute
and to the organizers of the Convex Geometry semester for their hospitality
and support.

\newpage \section{Preliminaries about Gaussian measure.}

We start with the

\medskip\noindent {\bf  Proof of Proposition 3. }
We follow the outline given in [I-MK] in the context of Komatsu
inequality.   Put $g(x) = e^{\frac{x^2}{2}}
\int_x^{\infty} e^{- \frac{t^2}{2}} dt$ and
$g_+ (x) = \frac{4}{3x + (x^2 +8)^{\frac{1}{2}}}$.  It is easily checked that
$g' = xg -1$ and somewhat more tediously verified that
$g_+' \leq xg_+ -1$.  Moreover, (e.g.) a direct calculation shows that
$g(x) \leq \frac{1}{x}$ for $x > 0$.  By considering the function $h = g_+
-g$ and
its differential inequality $h' \leq xh -1$ one gets (by the same argument
as in [I-MK]) that
$h= g_+ -g \geq 0$ on $(0,\infty$,  hence on $[0,\infty$.  Since $\lim_{x
\rightarrow -1}
g_+ (x) = \infty$ whereas $g( -1)$ is finite, it follows that $g_+ (x) \geq
g(x)$
also for $x \in (-1, 0)$  (otherwise consider $x \in (-1, 0)$ for which $h$
attains
its minimum).

The estimate from below is shown in a similar way  (and,  anyway, it is not
new). \hfill $\Box$\\

\medskip\noindent {\bf Remark 4.}  As was mentioned in the introduction,  for
$x > 0$
the upper estimate of Proposition 3 is sharper than the well known
estimate of Komatsu who proved that for $x > 0$
$$\frac{2}{x + (x^2 +4)^{\frac{1}{2}}} \leq e^{\frac{x^2}{2}} \int_x^{\infty}
e^{- \frac{t^2}{2}}dt \leq \frac{2}{x + (x^2 +2)^{\frac{1}{2}}}$$
(see [I-MK]).  We give below the values of relative ``errors" (rounded to
two significant digits) given by the two upper
estimates for some values of $x$;  we also list,  for reference,
the errors of the lower estimate.  All of these were calculated using
$Mathematica$ and verified with $Maple$.  Our estimate is clearly the tightest
of the three and vastly superior to the other upper estimate.

\bigskip

\begin{center}
\begin{tabular}{|c|c|c|c|} \hline
$x$ & Our Upper & Komatsu's Upper & Komatsu's Lower
\\ \hline \hline
$0$ & $.13$  &  $.13$   & $-.20$
\\  \hline
$2$ & $.30 \cdot 10^{-2}$ & $.67 \cdot 10^{-1}$ & $-.17 \cdot 10^{-1}$
\\ \hline
$4$  & $.20 \cdot 10^{-3}$ & $.25 \cdot 10^{-1}$ & $-.25 \cdot 10^{-2}$
\\ \hline
$6$ &  $.27 \cdot 10^{-4}$ & $.13 \cdot 10^{-1}$ & $-.61 \cdot 10^{-3}$
\\ \hline
$8$ &  $.59 \cdot 10^{-5}$ & $.74 \cdot 10^{-2}$ & $-.21 \cdot 10^{-3}$
\\  \hline
$10$ &  $.17 \cdot 10^{-5}$ & $.48 \cdot 10^{-2}$ & $-.92 \cdot 10^{-4}$
\\  \hline
$20$ &  $.30 \cdot 10^{-7}$ & $.12 \cdot 10^{-2}$ & $-.61 \cdot 10^{-5}$
\\  \hline
$30$ & $.27 \cdot 10^{-8}$ & $.55 \cdot 10^{-3}$ & $-.12 \cdot 10^{-5}$
\\  \hline
$40$ & $.48 \cdot 10^{-9}$ & $.31 \cdot 10^{-3}$ & $-.39 \cdot 10^{-6}$
\\  \hline
$50$ & $.13 \cdot 10^{-9}$ & $.20 \cdot 10^{-3}$ & $-.16 \cdot 10^{-6}$
\\  \hline
\end{tabular}

\bigskip \noindent {\em Relative errors of estimates for the ``Gaussian tail"
for selected values of $x$.}
\end{center}

\medskip The next result is a fairly easy consequence of Proposition 3.

\medskip \noindent {\bf Proposition 5.} {\em Let, for} $x \in {\bf R}$,
$$f(x) = \frac{e^{- \frac{x^2}{2}}}{\int_x^{\infty} e^{- \frac{t^2}{2}}
dt}.$$
{\em Then}
\newline  (i) $f(x)$ \mbox{ {\em is\  an increasing convex function.}}
\newline (ii) $x-f(x)$ \mbox{ {\em is an increasing (to $0$ as $x
\rightarrow \infty$) function.}}

\medskip\noindent {\bf Proof.} (i)
We compute
$$f'(x) = (\frac{ e^{- \frac{x^2}{2}}}{\int_{x}^{\infty}e^{-
\frac{t^2}{2}}dt})^2
-x\frac{ e^{- \frac{x^2}{2}}}{\int_{x}^{\infty}e^{- \frac{t^2}{2}}dt}$$
Clearly $f' \geq 0$ if and only if
$$\frac{ e^{- \frac{x^2}{2}}}{\int_{x}^{\infty}e^{- \frac{t^2}{2}}dt}-x
\geq 0.$$
If $x \leq 0$, this inequality holds trivially;  if $x > 0$ the inequality
holds e.g.    by Proposition 3, as $\frac{4}{3x + (x^2 +8)^{\frac{1}{2}}} \leq
\frac{1}{x}$.

\noindent We next have
$$f''(x) = \frac{1}{(e^{\frac{x^2}{2}}
\int_{x}^{\infty}e^{- \frac{t^2}{2}}dt)^3}
((x^2-1)(e^{\frac{x^2}{2}}
\int_{x}^{\infty}e^{- \frac{t^2}{2}}dt)^2 -3xe^{\frac{x^2}{2}}
\int_{x}^{\infty}e^{- \frac{t^2}{2}}dt + 2).$$

\noindent Clearly $f''(x) \geq 0$ if and only if
$$(x^2-1)(e^{\frac{x^2}{2}}
\int_{x}^{\infty}e^{- \frac{t^2}{2}}dt)^2 -3xe^{\frac{x^2}{2}}
\int_{x}^{\infty}e^{- \frac{t^2}{2}}dt + 2 \geq 0.$$

\noindent We put $z = e^{\frac{x^2}{2}}
\int_{x}^{\infty}e^{- \frac{t^2}{2}}dt$ and consider the expression
above as a polynomial in
$z$  i.e.  $z^2 (x^2 -1) - 3zx +2$.
As the roots of this polynomial are
$$z_{1/2} = \frac{3x \pm (x^2 +8)^{\frac{1}{2}}}{2(x^2 -1)} = \frac{4}{3x
\mp (x^2 +8)^{\frac{1}{2}}},$$
$f'' \geq 0$ holds trivially for $- \infty < x < -1$, and holds for
$x > -1$ if
$$e^{\frac{x^2}{2}} \int_x^{\infty} e^{- \frac{t^2}{2}} dt \leq
\frac{4}{3x + (x^2 +8)^{\frac{1}{2}}}$$
which is true by Proposition 3.

\par
(ii) By the calculation from the part (i) $$(x-f(x))'= 1 -
(\frac{ e^{- \frac{x^2}{2}}}{\int_{x}^{\infty}e^{- \frac{t^2}{2}}dt})^2
+ x\frac{ e^{- \frac{x^2}{2}}}{\int_{x}^{\infty}e^{- \frac{t^2}{2}}dt}.$$
After putting $z = \frac{e^{- \frac{x^2}{2}}}
{\int_{x}^{\infty}e^{- \frac{t^2}{2}}dt}$, the assertion $(x-f(x))'\geq 0$
becomes
$$1 + xz -z^2 \geq 0.$$
As the roots of this polynomial are
$$z_{1/2} = \frac{x \pm (x^2 +4)^{\frac{1}{2}}}{2},$$
the inequality follows, as before, from Proposition 3. \hfill $\Box$\\

\section{Proof of Theorem 1}

The proof of Theorem 1 is achieved in several steps.  In the first step we use
Ehrhard's inequality [E] to
reduce the general case to the $2$-dimensional case.  In the second step,
based on
(a rather general) Lemma 7, we reduce the $2$-dimensional problem even
further to
a four-parameter family of ``extremal" sets.  The final step is based on a
careful analysis of
dependence of the measures of sets involved on these parameters and uses
(computational) Lemmas 8 and 9.

Let $K$ and $u$ be as in Theorem 1 and let $H_0$ be
the hyperplane through $0$ orthogonal to $u$.  Without loss of generality
we may assume
that $\| u \|_2 \leq 1$.  For $t \in {\bf R}$ put
$H_t = H_0 + t \cdot u$ and let $\varphi (t) = \mu_{n -1} (K \cap H_t )$ and
$\Phi (x) = \mu_1 (( - \infty ,x])$.

\noindent By Ehrhard's inequality [E], $\psi (t) = \Phi^{-1} ( \varphi
(t))$ is a
concave function.
Therefore it is enough to consider the case $n =2$ and,  in place of $K$, sets
$K_{\psi}\subseteq {\bf R}^2$ of the form
\begin{equation}
K_{\psi} = \{ (x,y) \in {\bf R}^2 : y \leq \psi (x) \},
\end{equation}

\noindent with $u = e_1$ and $H_0$ identical with the $y$-axis,  where
$\psi$ is a
concave,
$\overline{{\bf R}}$ - valued function.  We will use the convention
$\Phi(- \infty)=0,
\Phi(\infty)=1$.  It may also be sometimes convenient to specify the interval
$[A,B]=\{x : \psi (x) > - \infty\}$.  The assumptions about the centroid become
\begin{equation}
\int_{\bf R} (x-c) \Phi ( \psi (x)) d \mu_1 (x) = 0 = \int_a^b (x-c)
d \mu_1 (x)
\end{equation}
and the assertion becomes
\begin{equation}
\int_a^b \Phi ( \psi (x)) d \mu_1 (x) \geq \int_{\bf R} \Phi ( \psi (x)) d
\mu_1 (x)
\int_a^b d \mu_1 (x).
\end{equation}

\medskip\noindent {\bf Remark 6.}   With this reduction of the general
case to the $2$-dimensional case we can now give a
quick proof of Sidak's Lemma.  As was indicated earlier,
Sidak's Lemma follows by induction
from the ``symmetric" variant of Theorem 1,  i.e.
when $L$ is a $0$-symmetric layer  ($a =-b$, $b > 0$) and $K$ is a
$0$-symmetric set
(hence $c =0$).  After reduction to the
$2$-dimensional case, $\psi$ is a concave function that is symmetric about
the $y$-axis (hence decreasing away from the origin) and one
has to show that
$$\int_{-b}^b \Phi ( \psi (x))d \mu_1 (x) \geq \int_{- \infty}^{\infty}
\Phi ( \psi (x))
d \mu_1 (x) \int_{-b}^b d \mu_1 (x)$$
or equivalently
$$\frac{\int_{-b}^b \Phi ( \psi (x)) d \mu_1 (x)}{\int_{-b}^b d \mu_1 (x)} \geq
\int_{- \infty}^{\infty} \Phi ( \psi (x)) d \mu_1 (x).$$
The above inequality holds because on the left we are averaging the function
$\Phi ( \psi (x))$ over the set where
it is ``biggest'',  while on the right - over the entire real line.

Actually it is not even necessary to use Ehrhard's inequality for this proof
of Sidak's
Lemma.  What is really used is (a special case of) the Brunn-Minkowski
inequality for Gaussian measure  (this was pointed out to the authors by A.
Giannopoulos) and the fact that the Gaussian measure is a product measure.

\medskip Returning to the proof of Theorem 1,   we show next that it is
enough to prove
inequality (5) for ``extremal" $\psi$'s which turn
out to be linear functions.  The reduction to this extremal case holds not only
for Gaussian measure on ${\bf R}^2$ but for a much more general class of
measures on ${\bf R}^2$ and is based on Lemma 7 that follows.  It will be
convenient to introduce
the following notation:  if  $\psi : [a,b] \rightarrow \overline{{\bf R}}$,
let
$$C_{\psi} = \{ (x,y): a \leq x \leq b,  y \leq \psi (x)\}.$$
We then have
\medskip \newline \noindent {\bf Lemma 7.}  {\em Let} $\psi : [a,b] \rightarrow
\overline{{\bf R}}$ {\em be a concave
function not identically  equal to} $- \infty$ {\em and let}
$\nu$ {\em be a finite measure on} ${\bf R}^2$ {\em that is absolutely
continuous with respect to
the Lebesgue measure.  Then there exists a linear function} $\psi_0 (x) =
mx +h$ {\em such that}

\medskip \noindent(i) \quad  $\nu (C_{\psi}) = \nu (C_{\psi_0})$
\newline (ii) \quad  $\int_{C_{\psi}} x d \nu (x,y) = \int_{C_{\psi_0}} x d
\nu (x,y)$
\newline (iii) \quad  $\psi_ (a) \leq \psi_0 (a),\quad \psi_ (b) \leq
\psi_0 (b)$
\newline (iv) \quad  $\psi_0' (a) \leq \psi'(a),\quad \psi_0' (b) \geq
\psi' (b)$

\bigskip We postpone the rather elementary proof of Lemma 7 until section 4.

\medskip\noindent  For $\alpha < \beta$ let us denote
$$L(\alpha,\beta) =  \{ (x,y) \in {\bf R}^2 : \alpha \leq x \leq \beta\}.$$
In the notation of Lemma 7 the assertion of Theorem 1  (or (5)) then becomes
\begin{equation}
\nu (C_{\psi}) \geq \nu (K_{\psi}) \cdot \nu (L(a,b)).
\end{equation}
\noindent (Note that $C_{\psi} = K_{\psi} \cap L(a,b)$;  the reader is
advised to draw a picture at this point to follow the remainder of the
argument).
\noindent Let now $\psi_0 (x) = mx+h$  be given by Lemma 7.  By symmetry,
we may assume
that $m \geq 0$.  The plan now is to show that,  for some (ultimately
unbounded) interval
$[A,B] \supset [a,b]$ and
\begin{equation}
$$\[ \psi_1(x) = \left \{ \begin{array}{ll}
   mx+h   & \mbox{if $x\in [A,B]$}\\
  -\infty  & \mbox{if $x\notin[A,B]$}
\end{array}
\right. \] $$
\end{equation}
we have
\begin{equation}
\nu (C_{\psi_1}) = \nu (C_{\psi})
\end{equation}
\begin{equation}
\nu (K_{\psi_1}) \geq \nu (K_{\psi})
\end{equation}
\noindent while,  at the same time,  the $\nu$-centroids of $K_{\psi_1}$ and
$K_{\psi}$
lie on the same line  $x=c$, i.e.
\begin{equation}
\int_{K_{\psi_1}} (x-c)d \nu (x,y) =\int_{\{A \leq x \leq B, y \leq mx+h\}}
(x-c)d \nu (x,y)
=\int_{K_{\psi}} (x-c)d \nu (x,y) =0.
\end{equation}
\noindent It will then follow immediately that it is enough to prove (6)
with $\psi$ replaced
by $\psi_1$, as required for reduction to the ``linear" case.
\noindent
\par
Now (8) is a direct consequence of the assertion (i) of Lemma 7 and (7).
On the other hand,
it follows from the assertions (iii) and (iv) that  $\psi_0(x)=mx+h \geq
\psi(x)$  for
$x \notin [a,b]$;  in other words $K_{\psi_0} \backslash L(a,b) \supset K_
{\psi}
\backslash L(a,b)$.
In combination with (8) this would imply  (9), {\em if} we were able
to set
$[A,B]=[-\infty,\infty]$.  However,  since we also need to ensure the
centroid assumption (10),
we need to proceed more carefully.
\noindent Let $A_0 \leq a$  (resp. $B_0 \geq b$)  be such that
\begin{equation}
\int_{K_{\psi_0}\cap L(A_0,a)} (x-c)d \nu (x,y) =\int_{K_{\psi}\cap
L(-\infty,a)} (x-c)d \nu (x,y).
\end{equation}
\noindent (resp. $L(b,\infty)$ and $L(b,B_0)$  in place of $L(-\infty,a)$
and $L(A_0,a)$).
This is possible since  $c \in (a,b)$ and,  as we indicated earlier,
$\psi_0 \geq \psi$  on
$(-\infty,a)$  (resp. on $(b,\infty)$).  Since,  by (i) and (ii) of Lemma 7,
$$\int_{C_{\psi_0}} (x-c)d \nu (x,y) =\int_{C_{\psi}} (x-c)d \nu (x,y),$$
\noindent it follows that the centroid condition (10) is satisfied if we
set  $[A,B]=[A_0,B_0]$.
Additionally,  an elementary argument shows that (11) combined with $\psi_0
\geq \psi$ on
$(-\infty,a]$  implies
$$\nu(K_{\psi_0} \cap L(A_0,a)) \geq  \nu(K_{\psi} \cap L(-\infty,a)).$$
\noindent This is roughly because the set on the left is ``closer" to the
axis $x=c$ than the one on the right and so,  for the ``moment equality"
(11)  to hold,  the former must have
a ``bigger mass".   Similarly,  $\nu(K_{\psi_0} \cap L(b,B_0)) \geq
\nu(K_{\psi} \cap L(b,\infty))$,
hence the ``mass condition"  (9)  also holds with  $[A,B]=[A_0,B_0]$.  This
reduces the problem
to linear functions  (more precisely functions of type  (7));  to get the
full reduction
(i.e. to an unbounded interval  $[A,B]$)  we notice that we may
simultaneously (and,
for that matter,  continuously)  move  $A$  to the left and  $B$ to the
right  starting from
$A_0, B_0$  respectively so that the centroid condition  (10)  holds,
until $A$  ``hits"
$-\infty$ or $B$  ``hits" $+\infty$;  the mass condition (9) will be then
{\em a fortiori} satisfied.

Thus,  depending on $c$, $m$  and $h$,  we end up with one of two possible
configurations
$$R_1=R_1 (h,B) = \{ (x,y) \in {\bf R}^2 : - \infty < x \leq B, y \leq mx+h \}$$
$$R_2=R_2 (h,A) = \{ (x,y) \in {\bf R}^2 : A \leq x < \infty , y \leq mx+h \},$$
for which we have, for $i = 1$ or $i=2$  (whichever applicable),
\begin{equation}
\mu_2 (R_i \cap L(a,b)) = \mu_2 (C_{\psi}) = \mu_2 (K_{\psi} \cap L(a,b))
\end{equation}
\begin{equation}
\mu_2 (R_i) \geq \mu_2 (K_{\psi})
\end{equation}
\begin{equation}
\int_{R_i} (x-c)d \mu_2 =\int_{K_{\psi}} (x-c)d \mu_2 = 0
\end{equation}
The three conditions above are just a rephrasing of  (8)-(10)  for  $\nu =
\mu_2$;
in particular it is enough to prove Theorem 1 for the extreme configurations
$K = R_i,  i=1,2$  or,  equivalently,  to prove  (5)  for  $\psi = \psi_1$
with
$\psi_1$ given by  (7) and some unbounded interval $[A,B]$. This will be
the last
step of the proof of the Theorem.

Let us note here that even though for the configuration $R_1=R_1 (h,B)$
it is possible in principle to have the centroid condition (14) satisfied
also for  $B<b$, we do not have to consider that case as it would have been
``reduced" in the previous step.  On the other hand,  one always has
$A \leq a$ for configurations of type $R_2$  (at least for  $m \geq 0$,
which we assume all the time).  See also the remarks following
the statement of Lemma 9.

For  $K=R_1$,   (5)  may be  restated as
\begin{equation}
\frac{\int_{- \infty}^B \Phi (mx +h)d \mu_1 (x)}{\int_a^b \Phi (mx +h)
\frac{d \mu_1 (x)}{\mu_1 ((a,b))}} \leq 1,
\end{equation}
while for $K=R_2$
\begin{equation}
\frac{\int_A^{\infty} \Phi (mx +h)d \mu_1 (x)}{\int_{a}^{b} \Phi (mx +h)
\frac{d \mu_1 (x)}{\mu_1 ((a, b))}} \leq 1.
\end{equation}

\medskip \noindent Denote the left hand side of (15) by $F_1 (h,w)$;
and the left hand side of
(16) by $F_2 (h,w)$,  where $w=\mu_1([a,b])$ is the ``Gaussian weight" of
the interval $[a,b]$.
Note that for fixed $c$ and $m$,
$B$ (resp. $A$) depends on $h$ as given by (14) with $i=1$ (resp. $i=2$).
Also note that it perfectly makes sense to consider
$h=+\infty$, $w=0$, $b=+\infty$ or $a=-\infty$ if otherwise allowable.

To study the behavior of $F_1$ and $F_2 $
we need two more lemmas.

\medskip\noindent {\bf Lemma 8.}  {\em With $B = B(h)$ (resp.  $A = A(h)$)
defined by (14) we have}
$$\frac{dB}{dh} \geq 0,\qquad    \frac{dA}{dh} \geq 0 \, .$$

\medskip\noindent {\bf Lemma 9.} {\em With $B = B(h)$ defined by (14) we have}
$$\frac{\partial F_1}{\partial h} (h,w) \geq 0  \, .$$

\medskip \noindent {\bf Proof of Lemma 8.}   We give the proof for $B$
(hence  $R_1
(h,B)$);
$A$ and $R_2$ are treated in a similar way.  Showing that
$\frac{dB}{dh} \geq 0$ for fixed $m$ and $c$ is equivalent to showing that
$\frac{dc}{dh} \leq 0$ for fixed $B$ and $m$.  Note that the centroid of
$R_1 (h,B)$ is a ``weighted average" of the centroids of the half lines
$y = mx + \overline{h}$, $- \infty < \overline{h} \leq h$, $- \infty < x
\leq B$.
Therefore to
show that $\frac{dc}{dh} \leq 0$ it is enough to show that the
$x$-coordinates of the centroids
of the halflines move further away from the line $x =B$ as $h$ increases.
We make a (orthogonal) change
of variable such that the line $y = mx +h$ becomes horizontal.  Denote the new
variables by $(u,v)$.  Showing that the $x$-coordinates of the centroids of the
halflines move further away from the line $x =B$ as $h$ increases is
equivalent to showing
that the $u$-coordinate of the centroids of the half-lines move further away
from the
corresponding value $U = U(h)$ on the line corresponding to $x =B$.  This
means that one has to
show that
$$U - \frac{\int_{- \infty}^U te^{- \frac{t^2}{2}} dt}{(2 \pi)^{\frac{1}{2}}
\Phi (U)}$$
increases as $U$ increases, which holds by Proposition 5 (ii). \hfill $\Box$\\
\par The computational proof of the Lemma 9 is somewhat involved;
we postpone it until the next section.

With Lemmas 8 and 9 we can conclude the proof of the Theorem.
Let us start with several observations concerning the qualitative dependence
of the regions  $R_i$ on $c$ and $h$  (for fixed $m>0$;  $m$ does not
{\em qualitatively} affect that dependence as long as it is positive,
the case $m=0$ being trivial).  These observations are only partly used in
the proof,
but they do clarify the argument nevertheless.  First,  if  $c<0$
(the special role of $0$ follows from the fact that the origin is the
centroid of
the entire plane),  then only configurations of type $R_1$ appear.  As $h$
increases,
$B=B(h)$ increases (by Lemma 8) and, as $h \rightarrow +\infty$, $B$
approaches some
limit value $\tilde{B}$  (of which we may think as $B(\infty)$)  defined by
the equation
$$\int_{-\infty}^{\tilde{B}} (x-c)d \mu_2 = 0$$
(cf. (14)).  It can also be shown that  as $h \rightarrow -\infty$, $B(h)$
approaches $c$,
but that has no bearing on our argument;  we do use only the fact that,
for fixed $w$ and $c$
(hence $a, b$), the condition $B(h) \geq b$ (on which we insist,  see the
remark following (16))
is,  again by Lemma 8, satisfied for $h$ in some interval (of the type
$[h^*,+\infty]$ if $c<0$).
If $c=0$,  the picture is similar except that  $B(\infty) = \infty$.
Finally,  if $c>0$,  $B(h)$  is also increasing with $h$,  except
that it reaches the limit value  $B=+ \infty$ for some finite $h=
\tilde{h}$,  at which
point the configuration  $R_2$  ``kicks in",  the half-plane $R_1
(\tilde{h}, +\infty)$
coinciding with $R_2 (\tilde{h}, -\infty)$.  As $h$ varies from $\tilde{h}$
to $+\infty$,
$A(h)$ increases from $-\infty$ to some limit value $\tilde{A}=A(\infty)$
defined by
$\int_{\tilde{A}}^{\infty} (x-c)d \mu_2 = 0$,  the limit set $R_2 (\infty,
A(\infty))$
being the half-plane  $\{ (x,y) \in {\bf R}^2 : x \geq A(\infty) \}$.

We first treat
$R_1 (h,B)$ when $c \leq 0$.  By Lemma 9, $\frac{\partial F_1}{\partial h}
\geq 0$ for all $w$.  Hence we are done in this case if we show (15) for
the extremal
configuration when $h = + \infty$ and $B=\tilde{B}$.  But then
$$R_1 ( \infty , \tilde{B}) = \{ (x,y): x \leq \tilde{B} \}$$
and hence
$$F_1 (h,w) \leq F_1 (\infty,w) = \mu_1 (- \infty, \tilde{B}) \leq 1$$
for all $h,w$.

Next we consider $R_1$ when $c \geq 0$.  In this case  Lemma 9 reduces the
deliberation to
the extremal configuration with $h= \tilde{h}$  (and $B=+ \infty$).  Now,
as we indicated earlier,
$$R_1 (\tilde{h},+ \infty ) = \{ (x,y) \in {\bf R}^2 : y \leq mx +
\tilde{h} \}=R_2 (\tilde{h}, -\infty )$$
and so the inequality (15) will follow if we show (16) with
$A=-\infty$ and the same values of $c, m$.
Thus it remains to handle the case of $R_2$ i.e. we have to show that
$$ F_2(h,w) \leq 1$$
for all h, w or equivalently that
\begin{equation}
\frac{\mu_2 (R_2(h,A)\cap L(a,b))}{\mu_2 (L(a,b))}
\geq  \mu_2 (R_2(h,A))
\end{equation}
for all h, w.
Now let us fix $h$ and $w$ ($m$ is fixed throughout the argument)
and vary $A$ (hence $c$). The right hand side of (17) is clearly largest
if $A = -\infty$. Similarly the left hand side is smallest if $A = -\infty$;
this follows from the fact that, as $A$ is decreasing to $-\infty$,
$c$ also decreases and consequently $L$ moves to the left so that
$\mu_2(R_2(h,A) \cap L(a,b))$ decreases. So also in the case
of $R_2$ we reduced the argument to the extremal configuration with
$A = -\infty$ and $h = \tilde{h}$. It remains to show that
\begin{equation}
\frac{\int_{-\infty}^{\infty}\Phi (mx + h) d\mu_1(x)}
{\int_{a}^{b}\Phi (mx + h) \frac{d\mu_1(x)}{\mu_1((a,b))}}
\leq 1.
\end{equation}

\noindent Throughout the remainder of the proof we will occasionally relax
the assumption that
$c$ is the Gaussian centroid of  $(a,b)$.  We first treat the case $h \geq 0$.
Observe that in that case $\frac{\int_{-d}^{d}\Phi  (mx + h) d\mu_ 1(x)}
{\int_{-d}^{d}d\mu_1(x)} $ decreases as d increases for $d \geq 0$ (this is
seen by
computing the derivative with respect to $d$).
Therefore
$$\int_{-\infty}^{\infty}\Phi (mx + h) d\mu_1(x)
\leq \frac{\int_{-d}^{d}\Phi (mx + h) d\mu_1(x)}
{\int_{-d}^{d}d\mu_1(x)}$$
The above is just (18) for $a=-b$. It now formally follows
that (18) holds whenever $ \frac{a+b}{2} \geq 0$ (or $b\geq -a$):
just compare the average of $\Phi (mx + h) $ over [a,b]
with that over $[-|a|,|a|]$ and use the fact that $\Phi (mx + h)$
is increasing in x.  In particular,  if $c$ {\em is} the Gaussian centroid
of  $(a,b)$,
then,  as is easily seen,  $ \frac{a+b}{2} \geq c \geq 0$,  which settles
the case $h \geq 0$.
\par
It remains to handle the case  $h < 0$.
\par
Let $\Phi_0 = \mu_2( \{ (x,y): y\leq mx+h\})$
and $h_0 = \Phi^{-1} (\Phi_0)$ (i.e. $\Phi_0 =\mu_2(\{(x,y): y\leq 0 \cdot
x+h_0\}))$.
We need to show that
\begin{equation}
\int_{a}^{b}\Phi (mx + h) \frac{d\mu_1(x)}{\mu_1((a,b))}
\geq \Phi_0 = \Phi(h_0).
\end{equation}

Let $x_0 = \frac{h_0-h}{m}$ be the x-coordinate of
the point of intersection of the lines $y = h_0$ and
$y = mx + h$. If $a \geq x_0$,
then (19) holds trivially,  hence we only need to consider the case $a < x_0$.
We will show that (19) holds provided $ \frac{a+b}{2} \geq x_0$. In our
situation
(i.e. when $c$ is the Gaussian centroid of  $(a,b)$) this condition is
satisfied since
 $\frac{a+b}{2} \geq c \geq x_0$.
Similarly as in the case of $h \geq 0$, it is enough to consider the case
$ \frac{a+b}{2} = x_0$  or  $b-x_0=x_0-a$. To show inequality (19), it is then
enough to show
\begin{equation}
\int_{a}^{x_0} (\Phi_0-\Phi(mx+h)) d\mu_1(x)
\leq \int_{x_0}^{b}
(\Phi(mx+h)-\Phi_0) d\mu_1(x)
\end{equation}
or equivalently, by rotational invariance of the Gaussian measure, that
\begin{equation}
\begin{array}{c}
\int_{x_0}^{x_1}
(\Phi(mx+h)-\Phi_0) d\mu_1(x) \\
+ \int_{x_1} ^{x_2}
(\Phi(\frac{x_0-x}{m}+ h_0 +
\frac{(1+m^2)^{(1/2)}}{m}(x_0-a))
-\Phi_0) d\mu_1(x) \nonumber \\
\leq
\int_{x_0}^{b}
(\Phi(mx+h)-\Phi_0) d\mu_1(x) ,
\end{array}
\end{equation}
where $x_1 = x_0 + \frac{x_0-a}{(1+m^2)^{(1/2)}}$
and $ x_2 = x_0 + (1+m^2)^{(1/2)}(x_0-a)$
(see Figure 1).

\begin{figure}
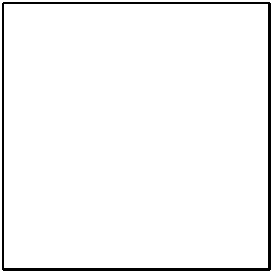

\BoxedEPSF{figuregauss.ps}
\caption{The case $h<0$.}
\end{figure}

\noindent Inequality (21) holds, if we can show that
$$\int_{b}^{x_2}
(\Phi(\frac{x_0-x}{m}+ h_0 +
\frac{(1+m^2)^{(1/2)}}{m}(x_0-a))-\Phi_0) d\mu_1(x)$$
$$\leq  \, \, \, \int_{x_1}^{b}
(\Phi(mx+h) -
\Phi(\frac{x_0-x}{m}+ h_0 +
\frac{(1+m^2)^{(1/2)}}{m}(x_0-a)) d\mu_1(x), $$

\noindent which holds as the triangles over which we integrate
have the same Lebesgue measure whereas the latter has bigger
Gaussian measure as the (restriction of the) reflection which maps the
first one into the
second is ``measure decreasing" with respect to the Gaussian measure.

As shown before, this also completes the proof of $F_1 \leq 1$ and
consequently that of the Theorem.  \hfill $\Box$\\

\medskip\noindent {\bf Remark 10.} We wish to reiterate that,  at least in
the case
when $K$ is a half plane  $\{(x,y) : y \leq mx+h\}$,  the requirement that
$c$ is the Gaussian
centroid of  $(a,b)$ may be relaxed somewhat:  to $ \frac{a+b}{2} \geq 0$
if  $h \geq 0$
and to $ \frac{a+b}{2} \geq x_0$ if $h \leq 0$. It follows that the same is true
for regions of type  $R_2$.  There is also some flexibility in the handling
of regions
of type $R_1$,  and consequently of an arbitrary $K$. However,  since we do
not have
any {\em natural} description of the allowed ``relaxation",  we do not
pursue this direction.

\section  {Proofs of the Lemmas.}

\medskip\noindent {\bf Proof of Lemma 7.}    We shall tacitly assume that
the density of $\nu$ with respect to the Lebesgue measure
is strictly positive, which is the case we need in our application;
the general case can be easily derived from this one.
We shall also assume that $\psi$ doesn't take the value $+\infty$ ,
in particular $\psi$  is continuous  (the opposite case is easy to handle
directly)
and that $\psi$ is not linear  (if it is,  we are already done).
For $m \in {\bf R}$ let the line $\psi^{(m)} (x) = mx +h$ be such that
\begin{equation}
\nu ( \{ (x,y): a \leq x \leq b, \, y \leq \psi (x) \} )
= \nu ( \{ (x,y): a \leq x \leq b, y \leq \psi^{(m)} (x) \} ), \nonumber
\end{equation}

\noindent where $h=h(m)$;  it follows from our assumptions that $h(\cdot )$
must be a continuous function.
The graph of $\psi^{(m)}$ cannot be completely above the graph of
$\psi$ on $(a,b)$ nor completely below the graph of $\psi$ on $(a,b)$;
otherwise the ``mass equality" (22) would not hold.
Therefore all the lines satisfying (22) {\em intersect} the graph of $\psi$
in at
least one point $(p, \psi (p))$ with $a < p < b$.

Now  suppose there is a line $\psi_0 (x) = m_0 x +h$ for which (22) holds,
for which the ``moment equality"
\begin{equation}
\int_{ \{ (x,y): a \leq x \leq b, y \leq \psi (x) \} } xd \nu =
\int_{ \{ (x,y): a \leq x \leq b, y \leq \psi_0 (x) \} } xd \nu
\end{equation}
holds and which has {\em exactly} one point of intersection $(p , \psi_0 (p))$
with the graph of $\psi$.   Then $\psi \leq \psi_0$ on one of the intervals
$[a,p],  [p,b]$
and $\psi \geq \psi_0$ on the other.  On the other hand,  it follows from
(22) and (23) that
$$\int_{ \{ (x,y): a \leq x \leq b, y \leq \psi (x) \} } (x -p) d \nu =
\int_{ \{ (x,y): a \leq x \leq b, y \leq \psi_0 (x) \} } (x -p) d \nu,$$

\noindent which is inconsistent with the preceding remark if $\psi$ and
$\psi_0$  are
not identical.  Consequently, the line $y= \psi_0(x)$ with the required
properties (22)
and (23) has to intersect the graph of
$\psi$ in at least two points $(p_1 , \psi (p_1))$, $(p_2, \psi (p_2))$ with
$a < p_1 < p_2 < b$ and, by concavity of $\psi$, in {\em exactly} two such
points.
Again by concavity of $\psi$ this is only possible if the assertions (iii)
and (iv)
of Lemma 7  hold.

It thus remains to show that among the linear functions $\psi^{(m)}$ for
which the
``mass equality" (22) (hence (i))  holds there is one for which also the
``moment equality"
(23) (hence (ii)) holds.  To this end,  observe that as $m \rightarrow +
\infty $,
the lines $y= \psi^{(m)}(x)$ ``converge" to a vertical line $x=a_1$,   where
$a_1$
is defined by $\nu (L(a_1, b)= \nu (C_{\psi})$.  One clearly has
$$\int_{C_{\psi}} xd \nu < \int_{L(a_1, b)} xd \nu.$$

Similarly,  as $m \rightarrow - \infty $,  the sets $C_{\psi ^{(m)}} $
``converge"
to a strip $L(a, b_1)$  satisfying $\int_{L(a, b_1)} xd \nu <
\int_{C_{\psi}} xd \nu$.
By continuity, there must be $m_0 \in {\bf R}$ such that $\psi_0=
\psi^{(m_0)}$ verifies (ii).
This finishes the proof of the Lemma.  \hfill $\Box$\\
\bigskip For the proof of Lemma 9 we shall need an elementary
auxiliary result.

\noindent {\bf Lemma 11.}  {\em Let} $g$ {\em be a convex function on an
interval} $[ \alpha , \beta ]$ {\em and let} $\rho$ {\em be a
positive measure on} $[\alpha , \beta ]$. {\em  Let } $\alpha',\beta'$
{\em be such that }
$\alpha \leq \alpha' < \beta' \leq \beta$ {\em and suppose that }
\begin{equation}
\frac{(\int_{\alpha}^{\beta} xd \rho (x)}{\rho ([\alpha, \beta])} -
\frac{\int_{\alpha'}^{\beta'}
x d \rho (x)}{\rho ([\alpha' , \beta'])})(g(\beta') - g(\alpha')) \geq 0.
\end{equation}
{\em Then}
$$\frac{\int_{\alpha}^{\beta} g (x) d \rho (x)}{\rho ([\alpha, \beta])}
\geq \frac{\int_{\alpha'}^{\beta'} g (x) d \rho (x)}{\rho ([\alpha' , \beta'])}.$$
\par
Note that if, in particular, $
\frac{(\int_{\alpha}^{\beta} xd \rho (x)}{\rho ([\alpha, \beta])} =
\frac{\int_{\alpha'}^{\beta'}
x d \rho (x)}{\rho ([\alpha' , \beta'])}$ or if
$g(\beta') = g(\alpha'))$, then the assertion holds.
We skip the proof (the reader is advised to draw a picture).

\medskip \noindent {\bf Proof of Lemma 9.}  We recall that by the comments
following the
statement of Lemma 9 (see also the remark preceding (15)), for
fixed $w$ and $c$
(hence fixed $a, b$), we do need to consider $h^* \leq h \leq \tilde{h}$,
where $h=h^*$ corresponds to
$B=b$ while  $\tilde{h}=+ \infty$ if $c \leq 0$ and $\tilde{h}$ ($<+
\infty$) is defined by
$B(\tilde{h})=+ \infty$ (or $A(\tilde{h})=- \infty$) if $c>0$.  We have
$$\frac{\partial F_1}{\partial h} = \frac{\int_{a}^{b} \Phi (mx +h)
\frac{d \mu_1 (x)}{\mu_1 (a,b)} \left[ B' e^{- \frac{1}{2} B^2} \Phi
(mB +h) + \int_{- \infty}^B e^{- \frac{1}{2} (mx +h)^2} d \mu_1 (x) \right]}
{(2 \pi)^{\frac{1}{2}} \left( \int_{a}^{b}
\Phi (mx +h) \frac{d \mu_1 (x)}{\mu_1 ((a,b))} \right)^2}$$

$$- \frac{\int_{- \infty}^B \Phi (mx +h) d \mu_1 (x) \int_{a}^{b}
e^{- \frac{1}{2} (mx +h)^2} \frac{d \mu_1 (x)}{\mu_1 (a,b)}}
{(2 \pi)^{\frac{1}{2}} \left( \int_{a}^{b} \Phi (mx +h) \frac{d \mu_1
(x)}{\mu_1 ((a,b))} \right)^2}.$$
As $B' \geq 0$ by Lemma 8,  $\frac{\partial F_1}{\partial h} \geq 0$ will
follow if
\begin{equation}
\frac{\int_{- \infty}^B \frac{e^{- \frac{1}{2} y^2}}{\Phi (y)} \Phi (y) d \mu_1
(x)}{\int_{- \infty}^B \Phi (y) d \mu_1 (x)} \geq \frac{\int_{a}^{b}
\frac{e^{- \frac{1}{2} y^2}}{\Phi (y)} \Phi (y) d \mu_1 (x)}{\int_{a}^{b}
\Phi (y) d \mu_1 (x)}
\end{equation}
where $y =mx +h$.

\noindent By Proposition 5, $g(y) = \frac{e^{- \frac{1}{2} y^2}}{\Phi (y)}$
is a convex decreasing function (note that $g(y) = f(-y)$,
where $f$ is as in Proposition 5. Moreover
$$\frac{\int_{- \infty}^B x\Phi (y) d \mu_1
(x)}{\int_{- \infty}^B \Phi (y) d \mu_1 (x)}=c=
\frac{\int_{a}^b x d \mu_1(x)}{\int_{a}^b d \mu_1 (x)}
\leq \frac{\int_{a}^{b}x\Phi (y) d \mu_1 (x)}{\int_{a}^{b}
\Phi (y) d \mu_1 (x)},$$
as $\Phi(y)$ is increasing,  and so the condition (24) is satisfied with
$d\rho(x) = \Phi (y) d \mu_1 (x)$. Consequently Lemma 11 yields (25),
completing the proof of Lemma 9.
\hfill $\Box$\\

\bigskip
\begin{tabular}{ll}
Stanislaw J. Szarek & Elisabeth Werner \\
Department of Mathematics & Department of Mathematics \\
Case Western Reserve University & Case Western Reserve University \\
Cleveland, Ohio 44106-7058 & Cleveland, Ohio 44106-7058 \\
E-mail: sjs13@po.cwru.edu & and \\
& Universit\'{e} de Lille \\
& UFR de Mathematiques \\
& Villeneuve d' Ascq, France  \\
& E-mail: emw2@po.cwru.edu \\
\end{tabular}


\begin{thebibliography}{M-PK.2} %
\bibitem[Ba]{ } R. Bagby: Calculating Normal Probabilities, Amer. Math.
Month. (1995), 46-49
\bibitem[B]{ } C. Borell: A Gaussian correlation inequality for certain bodies
in ${\bf R}\sp{n}$,  Math. Annalen 256 (1981), no. 4, 569-573.
\bibitem[E]{ } A. Ehrhard: Symetrisation dans l' espace de Gauss, Math.
Scand. 53 (1983),281-301.
\bibitem[G]{ } E.D. Gluskin: Extremal properties of orthogonal
parallelepipeds and their
application to the geometry of Banach spaces, Math. Sbornik vol. 64 (1985),
85-96.
\bibitem[H]{ } Y. Hu: Note on correlation and covariance inequalities, preprint.
\bibitem[I-MK]{ } K. Ito, H. P. McKean: Diffusion processes and their sample
paths,
Springer-Verlag, 1965.
\bibitem[K-MS]{ } A. L. Koldobsky, S. J. Montgomery-Smith: Inequalities of
correlation type for symmetric stable random vectors,
Statist. Probab. Lett. 28 (1996), No. 1, 91--97.
\bibitem[P]{ } L. Pitt: A Gaussian correlation inequality for symmetric
convex sets,
Annals of Probability 1977, vol. 5, No 3, 470-474.
\bibitem[S-S-Z]{ } G. Schechtman, T. Schlumprecht, J. Zinn: On the Gaussian
measure of the
intersection of symmetric convex sets, preprint.
\bibitem[S]{ } Z. Sidak: Rectangular confidence regions for the means of
multivariate normal distributions, J. Amer. Stat. Assoc. 62 (1967), 626-633.
%
\end{thebibliography}
\end{document}